\documentclass{amsart}
\usepackage{amssymb,latexsym}

\usepackage[T1]{fontenc}
\usepackage[utf8]{inputenc}
\usepackage{hyperref}
\usepackage{cleveref}
\usepackage{diagbox}
\usepackage{upgreek}
\usepackage{color}
\usepackage{bbm}
\usepackage{mathrsfs}
\usepackage{calligra}

\newtheorem{cor}{Corollary}

\newtheorem{lem}{Lemma}
\newtheorem{thm}{Theorem}

\newcommand{\bs}{\boldsymbol}

\newcommand{\Ausn}{\mathcal{P}_0}
\newcommand{\K}{\kappa}
\newcommand{\ka}{{\kappa+1}}

\allowdisplaybreaks

\title{Explicit bounds for prime $k$-tuplets}
\date{\today}
\author{Thomas Dubbe}

\begin{document}

\baselineskip=17pt

\begin{abstract}
Let $K\geq 2$ be a natural number and $a_i,b_i\in\mathbb{Z}$ for $i=1,\ldots,K-1$. We use the large sieve to derive explicit upper bounds for the number of prime $k$-tuplets, i.e., for the number of primes $p\leq x$ for which all $a_ip+b_i$ are also prime.
\end{abstract}

\maketitle

\section{Introduction}\label{S1}

For $K\in\mathbb{N}$ and $\bs{a},\bs{b}\in\mathbb{Z}^{K-1}$, we define $\pi_K(\bs{a},\bs{b};x)$ to be the number of primes $p\leq x$ such that $a_ip+b_i$ is also prime for every $1\leq i < K$. Famously, it is conjectured that we have $\pi_K \sim C_K x/\log^Kx$ for $x\rightarrow\infty$ and with a suitable constant $C_k$. Although the answer to this question, let alone the answer to the question of weather there is an infinite amount of prime $K$-tuplets, is still out of reach, we are able to obtain good upper estimates. In fact, we will show the following
\begin{thm}\label{THM1}
Let $K\geq 2$ be an integer and $\bs{a},\bs{b}\in\mathbb{Z}^{K-1}$ such that
\[ E_K = \prod_{i=1}^{K-1}a_ib_i\prod_{i<j}(a_ib_j-a_jb_i)\neq 0.\]
Furthermore let $\rho_K(p)$ be the number of solutions $\pmod{p}$ of the congruence
\[m\prod_{i=1}^{K-1}(a_im+b_i) \equiv 0\pmod{p}.\]
Also we require $\rho_K(p)<p$ for all $p\in\mathbb{P}$. Then we have 
\begin{align*}
\pi_K(\bs{a},\bs{b};X) &<2^K K!\prod_p\left(1-\frac{\rho_K(p)}{p}\right)\left(1-\frac{1}{p}\right)^{-K}\frac{X}{\log^KX}\left(\frac{\log X}{\log X - 2r_K}\right)
\end{align*}
for $X>\exp(2r_K)$ and
\[\pi_3(\bs{a},\bs{b};X) \leq 48\prod_p\left(1-\frac{\rho_3(p)}{p}\right)\left(1-\frac{1}{p}\right)^{-3}\frac{X}{\log^3X}\left(1 + \frac{34.549\log^2X}{P(\log X)}\right)\]
for $K=3$ and $X\geq 437400$, where $r_K$ and $P(t)$ are the same as in Theorem \ref{THM2}.
\end{thm}
Using Selberg's sieve one can show similar bounds for $\pi_K$ (see \cite{HALaRIC} and \cite{GRE}). However, these results are not fully explicit. The proof of our theorem generalizes an argument of Siebert who found upper bounds of this form for the special case $K=2$. His and our method of choice will be the large sieve. Therefore, we begin by repeating some basic notation. As usual, we have a set of natural numbers $\mathcal{A}\subset \mathbb{N}\cap [1,N]$ and a collection of sets $\Omega_p\subset\mathbb{Z}/p\mathbb{Z}$ consisting of residue classes $\pmod{p}$. For each $p$, the number of such classes is given by $\omega(p) = \vert\Omega_p\vert$. Then the sifting function is defined as  
\begin{align*} 
S(\mathcal{A},\Omega,z) = \#\{a\in\mathcal{A}\vert\; a\pmod{p} \notin \Omega_p\;\text{for all}\; p\leq z\}
\end{align*}
for some $z>0$. The set on the right is called \textit{sifted set}. Using the large sieve we can estimate its size through
\begin{lem}\label{Lem1}
By using the same notation as before we have
\begin{align*}
S(\mathcal{A},\Omega,z) \leq \frac{N}{G(z,\omega)}\quad\text{and}\quad S(\mathcal{A},\Omega,z) \leq \frac{N+z^2}{L(z,\omega)}.
\end{align*}
The denominators are given by 
\begin{align*}
G(z,\omega)&=\sum_{q\leq z}\mu^2(q)\left(1+\frac{3}{2}qzN^{-1}\right)^{-1}\prod_{p\mid q}\frac{\omega(p)}{p-\omega(p)},\\
L(z,\omega)&=\sum_{q\leq z}\mu^2(q)\prod_{p\mid q}\frac{\omega(p)}{p-\omega(p)}.
\end{align*}
Here $N>0$ is an integer and $z>0$ is an arbitrary positive number.
\end{lem}
The second inequality is the usual large sieve, while the first is its weighted version form Montgomery and Vaughan \cite{MONaVAU}. One gets better asymptotic estimates if one uses the latter at the expense of having weaker error terms.\\
We see that the essential information about the size of the sifted set is encoded in the numbers $\omega(p)$. So we have to make some assumptions regarding this function in order to proceed. For our purposes, we introduce the \textit{exception-set} $\Ausn$ as a finite subset of $\mathbb{P}$ and we require
\begin{align}\label{VOR1}
\omega(p) &< \kappa\quad \forall p\in \Ausn\notag\\
\omega(p) &= \kappa\quad \forall p\in \overline{\Ausn}=\mathbb{P}\setminus \Ausn, \tag{C1}
\end{align}
where $\kappa\geq 3$ is a natural number. We call $\K$ \textit{the dimension of the sieve}. In general, the dimension of a sieve indicates that $\omega(p)$ is on average equal to $\kappa$ (see again \cite{HALaRIC} or also \cite{FRIaIWA}). Moreover, we want to avoid the trivial case where every residue class is sifted out. So we demand
\begin{align}\label{VOR2}
\omega(p) < p\quad\forall p\in\mathbb{P}. \tag{C2}
\end{align}
Under this conditions we have
\begin{thm}\label{THM2}
Assume that (\ref{VOR1}) and (\ref{VOR2}) are true. Then we have on the one hand 
\begin{align*}
S(\mathcal{A},\Omega,z) <& 48\prod_p\left(1-\frac{\omega(p)}{p}\right)\left(1-\frac{1}{p}\right)^{-3}\frac{N}{P(\log N)}
\end{align*}
for $\kappa=3$, $N\geq 437 400$, $z=\sqrt{\frac{2}{3}N}$ and with  $P(t) = t^3 - 34.549t^2 + 379.423t -553.222$. On the other hand, for general $\K\geq 3$ we have the bound
\begin{align*}
S(\mathcal{A},\Omega,z) <2^{\K}\K!\prod_p\left(1-\frac{\omega(p)}{p}\right)\left(1-\frac{1}{p}\right)^{-\K}\frac{N}{\log^\K N}\left(1-\frac{r_\K}{\log N}\right)^{-1}
\end{align*}
for $N>\exp(r_\K)$, $r_\K = 5.961\cdot 10^{14}\K^2\log\K$ and with $z=\sqrt{\frac 2 3 N}$.
\end{thm}
We will prove this Theorem in the next three chapters using the weighted version of the large sieve from Lemma \ref{Lem1}. To do this, we look for lower estimates for $G(z,\omega)$, which we can derive by first treating $L(z,\omega)$ and then using summation by parts. Finally, we derive Theorem \ref{THM1} as an application of Theorem \ref{THM2}.


\section{The Term \texorpdfstring{$L(z,\omega)$}{L(z,w)}}\label{S2}

Our first task is to handle the effects of the exception set $\Ausn$. Instead of treating $L(z,\omega)$, we would prefer to deal with $L(z,\omega_0)$, where $\omega_0$ is as independent of $\Ausn$ as possible. Fortunately, we can transform the problem in such a way by the following

\begin{lem}\label{Lem2}
Let $f$ and $g$ be two functions with
\begin{enumerate}
\item[(i)] $p>f(p)\geq g(p)\geq 0\quad\forall p\in M$
\item[(ii)] $p>g(p)\geq f(p)\geq 0\quad\forall p\in\mathbb{P}\setminus M$.
\end{enumerate}
Here, $M$ is a finite set of primes. Then we have for all $z>0$
\begin{align*}
L(z,g) \geq \prod_{p\in M}\frac{p-f(p)}{p-g(p)}L(z,f).
\end{align*}
\end{lem}
The application now works as follows: we define the desired function $\omega_0$ by
\begin{align*}
\omega_0(p) = \begin{cases} \kappa & p>\kappa \\ \omega(p) & p\leq \kappa \end{cases}.
\end{align*}
Note that we have $\{p|p\leq \K\}\subseteq\Ausn$ because of (\ref{VOR1}) and (\ref{VOR2}). So we get after a quick examination that
\begin{align*}
p &> \omega_0(p)\geq \omega(p)\quad\forall p\in \Ausn \\
p &> \omega(p)=\omega_0(p) \quad\forall p\in\overline{\Ausn}.
\end{align*}
We deduce 
\begin{cor}\label{Cor1}
We have for $z>0$
\begin{align*}
L(z,\omega) \geq \prod_{\substack{p\in \Ausn\\ p>\kappa}}\frac{p-\kappa}{p-\omega(p)}L(z,\omega_0).
\end{align*}
\end{cor}
Proving the lemma is not too difficult.
\begin{proof}[Lemma 1]
We start by considering an arbitrary function $h$ with $h(p)<p$ for all $p\in\mathbb{P}$. We can express $L(h,y)$, for such $h$ as
\begin{align*}
L(y,h) =& \sum_{q\leq y}\mu^2(q)\prod_{p\mid q}\frac{h(p)}{p-h(p)} = \sum_{q\leq y}\mu^2(q)\prod_{p\mid q}\sum_{\nu=1}^\infty \left(\frac{h(p)}{p}\right)^{\nu}\\
       =& \sum_{q\leq y}\mu^2(q)\sum_{k(l)=q}\frac{1}{l}\prod_{p^\nu|| l}(h(p))^{\nu} = \sum_{k(l)\leq y} \frac{1}{l}\prod_{p^\nu|| l}(h(p))^{\nu}.
\end{align*}
Next, we define the arithmetic function $F$ by
\begin{align*}
F(p^{\nu})=\begin{cases} 1 & \nu=0\\ f(p)^{\nu-1}\cdot|f(p)-g(p)| & \nu\geq 1\end{cases}
\end{align*}
and the set $\mathcal{M}$ through
\begin{align*}
\mathcal{M} = \{n\in\mathbb{N}\mid n=\prod_{i\leq l}p_i^{\nu_i},\;p_i\in M,\nu_i\geq 0, l\in\mathbb{N}\}.
\end{align*}
This allows us rewrite the constant $A_M$ as 
\begin{align*}
A_M = \sum_{m\in\mathcal{M}}\frac{1}{m}F(m) = \sum_{m\in\mathcal{M}}\prod_{p^{\nu}||m}\frac{F(p^{\nu})}{p^{\nu}} = \prod_{p\in M}\sum_{\nu=0}^\infty \frac{F(p^{\nu})}{p^{\nu}} = \prod_{p\in M}\frac{p-g(p)}{p-f(p)}.
\end{align*}
With $\tilde{g}(l)=\prod_{p^{\nu}||l}(g(p))^\nu$ we deduce that
\begin{align*}
L(z,g)\cdot A_M & = \sum_{k(l)\leq z}\frac{1}{l}\tilde{g}(l)\cdot \sum_{m\in\mathcal{M}}\frac{1}{m}F(m) = \sum_{s}\frac{1}{s}\sum_{\substack{lm=s\\ k(l)\leq z\\ m\in\mathcal{M}}}\tilde{g}(l)F(m) \\
                 & \geq \sum_{k(s)\leq z}\frac{1}{s}\sum_{\substack{l\mid s\\ \frac{s}{l}\in\mathcal{M}}}\tilde{g}(l)F\left(\frac{s}{l}\right) = \sum_{k(s)\leq z}\frac{1}{s} G_M(s).
\end{align*}
It can easily be shown that the function
\[
G_M(s) = \sum_{\substack{l\mid s\\ \frac{s}{l}\in\mathcal{M}}}\tilde{g}(l)F\left(\frac{s}{l}\right)
\]
is multiplicative and fulfils
\[
G_M(p^\nu) \geq (f(p))^\nu, \quad\forall p\in\mathbb{P}.
\]
This finally yields
\begin{align*}
L(z,g)\cdot A_M \geq \sum_{k(s)\leq z}\frac{1}{s}G_M(l) \geq \sum_{k(l)\leq z}\frac{1}{l} \prod_{p^\nu|| l}(f(p))^\nu = L(z,f),
\end{align*}
but this is equivalent to our original statement.
\end{proof}


\section{Reduction to and estimation of \texorpdfstring{$F_{\kappa}(z)$}{Fk(z)}}\label{S3}
\subsection{The main lemma and its proof}
Although $\omega_0$ is easier to work with than our original $\omega$, the behaviour for $p\leq \kappa$ still remains a problem. Therefore, we separate those special cases by writing
\begin{align}\label{EQ1}
L(z,\omega_0) &= \sum_{q\leq z}\mu^2(q) \prod_{p\mid q}\frac{\omega_0(p)}{p-\omega_0(p)} = \sum_{q\leq z}\mu^2(q) \prod_{\substack{p\mid q\\ p>\kappa}}\frac{\kappa}{p-\kappa}\prod_{\substack{p\mid q\\ p\leq \kappa}}\frac{\omega(p)}{p-\omega(p)} \notag\\
&= \sum_{l\mid P(\kappa)}\prod_{p\mid l}\frac{\omega(p)}{p-\omega(p)}F_{\kappa}\left(\frac{z}{l}\right). 
\end{align}
Here the function $F_k(y)$ is defined as
\begin{align*}
F_k(y) = \sum_{\substack{q\leq y\\ (q,P(k))=1}}\mu^2(q)\prod_{p\mid q}\frac{k}{p-k}.
\end{align*}
Since $F_k(y)$ is obviously monotonically increasing in $y$, we get for (\ref{EQ1}) 
\begin{align}\label{EQ2}
L(z,\omega_0) \geq F_{\kappa}\left(\frac{z}{P(\kappa)}\right)\sum_{l\mid P(\kappa)}\prod_{p\mid l}\frac{\omega(p)}{p-\omega(p)} = \prod_{p\leq \kappa}\frac{p}{p-\omega(p)} F_{\kappa}\left(\frac{z}{P(\kappa)}\right).
\end{align}
Once again, we have simplified our problem. This time we need lower bounds for $F_\kappa(z)$ in order to estimate $L(z,\omega_0)$. We can achieve this goal through the following lemma
\begin{lem}\label{Lem3}
We have for $z\geq 1$ 
\begin{align*}
F_3(z)\geq W_3\left(\frac{1}{6}\log^3 z + C_{1,3} \log^2 z + C_{2,3} \log z + C_{3,3}\right),
\end{align*}
with $C_{1,3} = -1.5351$, $C_{2,3} = 4.1864$ and $C_{3,3} = 16.7665$. For $\K\geq 3$ and $z\geq 1$ we get 
\begin{align*}
F_{\K}(z)\geq W_\K\left(\frac{1}{\K!}\log^\K z - C_\K\log^{\K-1}z\right).
\end{align*}
The constants are given by $C_\K = \frac{\exp(2e)}{(\K-1)!}\sum_{\ell=3}^{\K}\exp\left(\frac{e^2}{\ell}\left(\frac{\log\ell}{\log{\ell}-1}\right)\right)\log\ell$ and
\[ W_\K = \prod_p \left(1-\frac{\chi_\K(p)}{p}\right)^{-1}\left(1-\frac{1}{p}\right)^{\K}\prod_{p\leq\K}\left(1-\frac{1}{p}\right),\]
with $\chi_\K(p)=1$ for $p\leq\K$ and $\chi_\K(p) = \K$ otherwise.
\end{lem}

Before beginning the proof, we turn our attention to a lemma that will be quite useful for our purposes.

\begin{lem}\label{Lem4}
Let $f(n)\geq 0$ be a multiplicative function and $k,k_0$ be co prime integers. Then we have for $x\geq 1$
\begin{align*}
\sum_{\substack{l\leq x \\ (l,k_0)=1 \\ (l,k)=1}}\mu^2(l)f(l) \geq \prod_{p\mid k}\left(1+f(p)\right)^{-1}\sum_{\substack{l\leq x \\ (l,k_0)=1}} \mu^2(l)f(l).
\end{align*}
\end{lem}

This result is quite well known. Nevertheless, we want to give a short proof for the sake of completeness
\begin{proof}
Let $\epsilon_k(n)$ be $1$ for $(n,k)=1$ and $0$ otherwise. Furthermore, we define $g(n) = f(n)\epsilon_{k_0}(n)$ and
\[
G_k(x) = \sum_{\substack{l\leq x \\ (l,k)=1}}\mu^2(l)g(l),\quad G(x)=G_1(x).
\]
It follows that
\begin{align*}
G(x) & = \sum_{l\mid k}\sum_{\substack{m\leq x \\ (m,k)=l}}\mu^2(m)g(m) = \sum_{l\mid k}\sum_{\substack{n\leq x/l \\ (n,k/l)=1 \\ (n,l)=1}}\mu^2(l\cdot n)g(l\cdot n) \\
     & = \sum_{l\mid k}\mu^2(l)g(l)G_k\left(\frac{x}{l}\right)\leq\left(\sum_{l\mid k}\mu^2(l)g(l)\right)G_k(x).
\end{align*}
Moreover we have
\begin{align*}
\sum_{l\mid k}\mu^2(l)g(l) = \prod_{l\mid k}(1+g(p)) \leq \prod_{l\mid k}(1+f(p)),
\end{align*}
which gives us
\begin{align*}
\sum_{\substack{l\leq x \\ (l,k)=1}}\mu^2(l)g(l) \geq \prod_{l\mid k}(1+f(p))^{-1}\sum_{l\leq x}\mu^2(l)g(l).
\end{align*}
The statement follows by inserting $g(n) = f(n)\epsilon_{k_0}(n)$.
\end{proof}

\begin{proof}
First, we want to prove the recursive relation
\[ F_{\ka}(z) \geq \frac{P(\K)}{P(\ka)}\prod_{p\leq \ka}\left(1-\frac{1}{p}\right)\sum_{\substack{\ell\leq z\\(l,P(\ka))=1}}\mu^2(\ell)\frac{\tau_\K(\ell)}{\psi_\ka(\ell)\ell}\int_1^{z/\ell}\frac{F_\K(t)}{t}\;dt\]
for arbitrary $\K\geq 1$. This handy estimate will help us to derive our Lemma inductively. We begin by writing 
\begin{align*}
F_{\K+1}(z) = \sum_{\substack{d\leq z\\ (d,P(\K+1))=1}}\mu^2(d)\prod_{p\mid d}\frac{\K+1}{p-(\K+1)} = \sum_{\substack{d\leq z\\ (d,P(\K+1))=1}}\mu^2(d)\frac{\tau_{\K+1}(d)}{\psi_{\K+1}(d)}.
\end{align*}
Here $\tau_n(m) = \#\{(d_1,\ldots,d_n)\in\mathbb{N}^n\mid d_1\cdot\ldots\cdot d_n = m\}$ is the generalized divisor function and $\psi_n(m) = \prod_{\substack{p\mid m\\ p>n}} (p-n)$ is multiplicative. The same can be said about the former since  $\tau(m) = \sum_{d\mid m}1$ is multiplicative and we have the recursion  $\tau_n(m) = \sum_{d\mid m}\tau_{n-1}(d)$. Repeated use of lemma \ref{Lem4} brings us to
\begin{align*}
F_\ka(z) &= \sum_{\substack{d\leq z\\ (d,P(\ka))=1}} \mu^2(d)\frac{\tau_\ka(d)}{\psi_\ka(d)} = \sum_{\substack{d\leq z\\ (d,P(\ka))=1}} \frac{\mu^2(d)}{\psi_\ka(d)}\sum_{m\mid d}\tau_\K\left(\frac{d}{m}\right)\\
         =& \sum_{\substack{d\leq z\\ (d,P(\ka))=1}} \frac{\mu^2(d)}{\psi_\ka(d)}\sum_{\substack{\ell\leq z/d\\ (\ell,P(\ka))=1 \\ (d,\ell)=1}}\mu^2(\ell)\frac{\tau_\K(\ell)}{\psi_\ka(\ell)}\\
         \geq& \sum_{\substack{d\leq z\\ (d,P(\ka))=1}} \frac{\mu^2(d)}{\psi_\ka(d)} \prod_{p\mid d}\left(1+\frac{\tau_\K(p)}{\psi_\ka(p)}\right)^{-1}\sum_{\substack{\ell\leq z/d\\ (\ell,P(\ka))=1}}\mu^2(\ell)\frac{\tau_\K(\ell)}{\psi_\ka(\ell)}\\
         =& \sum_{\substack{d\leq z\\ (d,P(\ka))=1}} \frac{\mu^2(d)}{\varphi(d)}\sum_{\substack{\ell\leq z/d\\ (\ell,P(\ka))=1}}\mu^2(\ell)\frac{\tau_\K(\ell)}{\psi_\K(\ell)}\cdot\frac{\psi_\K(\ell)}{\psi_\ka(\ell)}\\
         =& \sum_{\substack{d\leq z\\ (d,P(\ka))=1}} \frac{\mu^2(d)}{\varphi(d)}\sum_{\substack{\ell\leq z/d\\ (\ell,P(\ka))=1}}\mu^2(\ell)\frac{\tau_\K(\ell)}{\psi_\K(\ell)}\prod_{p\mid \ell}\left(1+\frac{1}{\psi_\ka(p)}\right)\\
         =& \sum_{\substack{d\leq z\\ (d,P(\ka))=1}} \frac{\mu^2(d)}{\varphi(d)}\sum_{\substack{\ell\leq z/d\\ (\ell,P(\ka))=1}}\mu^2(\ell)\frac{\tau_\K(\ell)}{\psi_\K(\ell)}\sum_{m\mid \ell}\frac{\mu^2(m)}{\psi_\ka(m)}\\
         =& \sum_{\substack{d\leq z\\ (d,P(\ka))=1}} \frac{\mu^2(d)}{\varphi(d)}\sum_{\substack{\ell\leq z/d\\ (\ell,P(\ka))=1}}\mu^2(\ell)\frac{\tau_\K(\ell)}{\psi_\ka(\ell)\psi_\K(\ell)}\times\\
         &\quad\sum_{\substack{n\leq z/d\ell\\ (n,P(\ka))=1\\ (n,\ell)=1}}\mu^2(n)\frac{\tau_\K(n)}{\psi_\K(n)}\\
         \geq& \sum_{\substack{d\leq z\\ (d,P(\ka))=1}} \frac{\mu^2(d)}{\varphi(d)}\sum_{\substack{\ell\leq z/d\\ (\ell,P(\ka))=1}}\mu^2(\ell)\frac{\tau_\K(\ell)}{\psi_\ka(\ell)\ell}\sum_{\substack{n\leq z/d\ell\\ (n,P(\ka))=1}}\mu^2(n)\frac{\tau_\K(n)}{\psi_\K(n)}\\
         \geq& \frac{P(\K)}{P(\ka)}\sum_{\substack{d\leq z\\ (d,P(\ka))=1}} \frac{\mu^2(d)}{\varphi(d)}\sum_{\substack{\ell\leq z/d\\ (\ell,P(\ka))=1}}\mu^2(\ell)\frac{\tau_\K(\ell)}{\psi_\ka(\ell)\ell}F_\K\left(\frac{z}{d\ell}\right).
\end{align*}
Changing the order of summation, estimating $\sum_{d\leq x}\frac{\mu^2(d)}{\varphi(d)} \geq \log x$ (valid for $x \geq 1$, see \cite{MONaVAU}) and using lemma \ref{Lem4}, we get the desired result
\begin{align*}
&= \frac{P(\K)}{P(\ka)}\sum_{\substack{d\leq z\\ (d,P(\ka))=1}} \frac{\mu^2(d)}{\varphi(d)}\sum_{\substack{\ell\leq z/d\\ (\ell,P(\ka))=1}}\mu^2(\ell)\frac{\tau_\K(\ell)}{\psi_\ka(\ell)\ell}\sum_{\substack{n\leq z/d\ell \\ (n,P(\K))=1}}\mu^2(n)\frac{\tau_\K(n)}{\psi_\K(n)}\\
         &= \frac{P(\K)}{P(\ka)}\sum_{\substack{\ell\leq z\\ (\ell,P(\ka))=1}}\mu^2(\ell)\frac{\tau_\K(\ell)}{\psi_\ka(\ell)\ell}\sum_{\substack{n\leq z/\ell \\ (n,P(\K))=1}}\mu^2(n)\frac{\tau_\K(n)}{\psi_\K(n)}\sum_{\substack{d\leq z/\ell n\\ (d,P(\ka))=1}} \frac{\mu^2(d)}{\varphi(d)}\\
         &\geq \frac{P(\K)}{P(\ka)}\prod_{p\leq \ka}\left(1-\frac{1}{p}\right)\sum_{\substack{\ell\leq z\\ (\ell,P(\ka))=1}}\mu^2(\ell)\frac{\tau_\K(\ell)}{\psi_\ka(\ell)\ell}\times\\
         &\quad\quad\sum_{\substack{n\leq z/\ell \\ (n,P(\K))=1}}\mu^2(n)\frac{\tau_\K(n)}{\psi_\K(n)}\sum_{d\leq z/\ell n} \frac{\mu^2(d)}{\varphi(d)}\\
         &\geq \frac{P(\K)}{P(\ka)}\prod_{p\leq \ka}\left(1-\frac{1}{p}\right)\sum_{\substack{\ell\leq z\\ (\ell,P(\ka))=1}}\mu^2(\ell)\frac{\tau_\K(\ell)}{\psi_\ka(\ell)\ell}\times\\
         &\quad\quad\sum_{\substack{n\leq z/\ell \\ (n,P(\K))=1}}\mu^2(n)\frac{\tau_\K(n)}{\psi_\K(n)}\log\frac{z}{\ell n}\\
         &= \frac{P(\K)}{P(\ka)}\prod_{p\leq \ka}\left(1-\frac{1}{p}\right)\sum_{\substack{\ell\leq z\\ (\ell,P(\ka))=1}}\mu^2(\ell)\frac{\tau_\K(\ell)}{\psi_\ka(\ell)\ell}\int_1^{z/\ell}\frac{F_\K(t)}{t} dt.
\end{align*}
We proceed as planned by proving our Lemma by induction over $\kappa\geq 3$. Therefore, we start with the case $\kappa = 3$. In fact, the inequality we will show is slightly stronger than needed for the induction. However, we take a little extra care in this case, because we get sharper estimates for the number of prime $3-$ tuple and it may be needed for some future applications. Our starting point is the recursion
\[ F_3(z) \geq \frac{1}{9}\sum_{\substack{\ell\leq z \\ (l,6)=1}}\mu^2(\ell)\frac{\tau(\ell)}{\psi_3(\ell)\ell}\int_1^{z/\ell}\frac{F_2(t)}{t} dt.\]
According to Siebert \cite{SIE} (page $332$, ineq. $(2.7)$) we have
\[ F_2(z) \geq \frac{1}{4\omega}\left(\frac{1}{2}\log^2 z + C_{1,2}\log z + C_{2,2}\right),\]
for $z\geq 1$ with $C_{1,2}= 1.291$ and $C_{2.2} = 1.169$ (Siebert calculated the slightly wrong value $C_{2,2}= 1.166$). Furthermore, $\omega = \prod_{p>2}\left(1-\frac{1}{(p-1)^2}\right) = 0.660161815...$ is the well known twin prime constant. We deduce for $F_3$ the following
\begin{align*}
F_3(z) &\geq \frac{1}{36\omega_2} \sum_{\substack{\ell\leq z \\ (l,6)=1}}\mu^2(\ell)\frac{\tau(\ell)}{\psi_3(\ell)\ell}\int_1^{z/\ell}\frac{1}{2}\log^2 t + C_{1,2}\log t + C_{2,2}\;\frac{dt}{t} \\
       &= \frac{1}{36\omega_2} \sum_{\substack{\ell\leq z \\ (l,6)=1}}\mu^2(\ell)\frac{\tau(\ell)}{\psi_3(\ell)\ell}\left\{\frac{1}{6}\log^2 \frac{z}{\ell} + \frac{C_{1,2}}{2}\log^2 \frac{z}{\ell} + C_{2,2}\log\frac{z}{\ell}\right\}.
\end{align*}
Next, we set 
\[S_3(z) = \sum_{\substack{d>z\\ (d,6)=1}}\mu^2(d)\frac{\tau(d)}{\psi_3(d)d} \]
and $\frac{1}{\alpha_3} = S_3(0)$. We also require that $z$ is greater than or equal to a parameter $z_0>1$, which we will choose later. It follows that
\begin{align*}
F_3(z) &\geq \frac{1}{36\omega_2} \sum_{\substack{\ell\leq z \\ (l,6)=1}}\mu^2(\ell)\frac{\tau(\ell)}{\psi_3(\ell)\ell}\left\{\frac{1}{6}\log^2 \frac{z}{\ell} + \frac{C_{1,2}}{2}\log^2 \frac{z}{\ell} + C_{2,2}\log\frac{z}{\ell}\right\}\\
       &= \frac{1}{36\omega_2}\int_1^z \left(\sum_{\substack{\ell\leq t\\ (\ell,6)=1}}\mu^2(\ell)\frac{\tau(\ell)}{\psi_3(\ell)\ell}\right)\left(\frac{1}{2}\log^2\frac{z}{t} + C_{1,2}\log\frac{z}{t} + C_{2,2}\right)\frac{dt}{t}\\
			 &\geq \frac{1}{36\omega_2} \sum_{\substack{\ell\leq z_0 \\ (l,6)=1}}\mu^2(\ell)\frac{\tau(\ell)}{\psi_3(\ell)\ell}\left\{\frac{1}{6}\log^2 \frac{z_0}{\ell} + \frac{C_{1,2}}{2}\log^2 \frac{z_0}{\ell} + C_{2,2}\log\frac{z_0}{\ell}\right\}\\
			 &\quad +\int_{z_0}^z \left(\frac{1}{\alpha_3}-S_3(t)\right)\left(\frac{1}{2}\log^2\frac{z}{t} + C_{1,2}\log\frac{z}{t} + C_{2,2}\right)\frac{dt}{t}\\
			 &=\frac{1}{36\omega_2\alpha_3}\Big(\frac{1}{6}\log^3\frac{z}{z_0} + \frac{C_{1,2}}{2}\log^2\frac{z}{z_0} + C_{2,2}\log\frac{z}{z_0} \\
			 &\quad + \alpha_3 \sum_{\substack{\ell\leq z_0 \\ (l,6)=1}}\mu^2(\ell)\frac{\tau(\ell)}{\psi_3(\ell)\ell}\left\{\frac{1}{6}\log^2 \frac{z_0}{\ell} + \frac{C_{1,2}}{2}\log^2 \frac{z_0}{\ell} + C_{2,2}\log\frac{z_0}{\ell}\right\}\\
			 &\quad - \alpha_3\int_{z_0}^z S_3(t)\left(\frac{1}{2}\log^2\frac{z}{t} + C_{1,2}\log\frac{z}{t} + C_{2,2}\right)\frac{dt}{t}\Big).
\end{align*}
We can estimate $S_3(t)$ by
\[ S_3(t) < \frac{1}{\alpha_3}\frac{1534}{42875}\frac{\log^2 t + 45.75 \log t + 524.266}{t}\]
for $t\geq 1$. We will show the inequality right after this proof in the next subsection. Using this, we obtain the following for $F_3(z)$
\begin{align*}
F_3(z) &\geq \frac{1}{36\omega_2\alpha_3}\Big(\frac{1}{6}\log^3\frac{z}{z_0} + \frac{C_{1,2}}{2}\log^2\frac{z}{z_0} + C_{2,2}\log\frac{z}{z_0} \\
       &\quad + C_0(z_0) - \frac{1534}{42875}I(z,z_0)\Big).
\end{align*}
The sum $C_0(z)$ is given by 
\[ C_0(z) = \alpha_3 \sum_{\substack{\ell\leq z \\ (l,6)=1}}\mu^2(\ell)\frac{\tau(\ell)}{\psi_3(\ell)\ell}\left\{\frac{1}{6}\log^2 \frac{z}{\ell} + \frac{C_{1,2}}{2}\log^2 \frac{z}{\ell} + C_{2,2}\log\frac{z}{\ell}\right\},\]
while the integral yields
\[ I(z,z_0) = \int_{z_0}^z \frac{\log^2 t + 45.75\log t + 524.266}{t}\cdot\frac{\frac{1}{2}\log^2\frac{z}{t} + C_{1,2}\log\frac{z}{t} + C_{2,2}}{t}dt.\]
The latter can be further evaluated to 
\begin{align*}
&I(z,z_0) \leq \frac{1}{z}\left(-0.875\log^2z-43.3425\log z - 539.502\right)\\
             & -\frac{1}{z_0}\bigg(\left(-\frac{1}{2}\log^2z_0 - 23.875\log z_0 - 286.008\right)\log^2 z\\
					   & +\left(\log^3 z_0 + 47.459\log^2 z_0 + 560.1207\log z_0 -116.707\right)\log z\\
					   & -\frac{1}{2}\log^4z_0 - 23.584\log^3z_0 - 274.9908\log^2 z_0 + 73.3644\log z_0 - 539.503\bigg).
\end{align*}
Finally, we can estimate $F_3$ by
\begin{align*}
F_3(z)\geq& \frac{1}{36\omega_2\alpha_3}\Big(\frac{1}{6}\log^3\frac{z}{z_0} + \frac{C_{1,2}}{2}\log^2\frac{z}{z_0} + C_{2,2}\log\frac{z}{z_0} \\
          & + C_0(z_0) - \frac{1534}{42875}I(z,z_0)\Big)\\
			\geq& W_3\bigg(\frac{1}{6}\log^3z + A_{1,3}(z_0)\log^2z + A_{2,3}(z_0)\log z + A_{3,3}(z_0)\\
			    &+\frac{1534}{42875}\frac{1}{z}\left(-0.875\log^2 z - 43.3425\log z - 539.503\right)\bigg)\\
			\geq& W_3\bigg(\frac{1}{6}\log^3z + A_{1,3}(z_0)\log^2z + A_{2,3}(z_0)\log z + A_{3,3}(z_0),
\end{align*}
for $z\geq z_0$, with the constants 
\begin{align*}
&A_{1,3}(z_0) = \frac{C_{1,2}}{2} - \frac{\log z_0}{2} - \frac{1534}{42875}\frac{1}{z_0}\left(\frac{1}{2}\log^2z_0 + 23.875\log z_0 + 286.008\right)\\
&A_{2,3}(z_0) = \frac{\log^2 z_0}{2}-C_{1,2}\log z_0 + C_{2,2}\\
              &\quad\quad+ \frac{1534}{42875}\frac{1}{z_0}\left(\log^3 z_0 + 47.459\log^2 z_0 + 560.1207\log z_0 - 116.707\right)\\
&A_{3,3}(z_0) = \frac{C_{1,2}}{2}\log^2 z_0 - C_{2,2}\log z_0 + C_0(z_0)- \frac{1534}{42875}\frac{1}{z_0}\times\\
              &\left(\frac{1}{2}\log^4 z_0 +23.584\log^3 z_0 + 274.9908\log^2 z_0 - 73.364 \log z_0 + 539.503\right).
\end{align*}
Choosing $z_0 = 39$ gives us the first part of our lemma for $z\geq 39$. The remaining cases can be verified by direct computation. \\
We continue to complete the induction. To do this, we assume that the statement holds for some arbitrary $\K\geq 3$ and we define the polynomials $P_k(t) = \frac{1}{k!}t^k - C_kt^{k-1}$ and $Q_k(t) = \int_0^t P_k(t')\;dt'$, as well as $C'_k = P(k-1)/P(k)\prod_{p\leq k}\left(1-1/p\right)$. Analogous to the case $\K = 3$ we also set 
\[ S_k(z) = \sum_{\substack{d>z\\ (d,P(k))=1}}\mu^2(d)\frac{\tau_{k-1}(d)}{\psi_k(d)d},\]
$\frac{1}{\alpha_k} = S_k(0)\geq 1$ and $S'_k(t) = \frac{1}{\alpha_k} - S_k(t)$. Using the recursive relation from the beginning and the induction hypothesis, we get
\begin{align}\label{EQ3}
F_\ka(z) &\geq \frac{P(\K)}{P(\ka)}\prod_{p\leq \ka}\left(1-\frac{1}{p}\right)\sum_{\substack{\ell\leq z\\(l,P(\ka))=1}}\mu^2(\ell)\frac{\tau_\K(\ell)}{\psi_\ka(\ell)\ell}\int_1^{z/\ell}\frac{F_\K(t)}{t}\;dt\notag\\
         &\geq C'_\ka W_\K\sum_{\substack{\ell\leq z\\ (\ell,P(\ka))=1}}\mu^2(\ell)\frac{\tau_\K(\ell)}{\psi_\ka(\ell)\ell}\int_1^{z/\ell}\frac{P_\K(\log(t))}{t} dt\notag\\
         &= C'_\ka W_\K\sum_{\substack{\ell\leq z\\ (\ell,P(\ka))=1}}\mu^2(\ell)\frac{\tau_\K(\ell)}{\psi_\ka(\ell)\ell}Q_\K\left(\log\frac{z}{\ell}\right)\notag\\
         &= C'_\ka W_\K \int_1^z S'_\ka(t)\frac{P_\K\left(\log\frac{z}{t}\right)}{t}\;dt\notag\\
         &= C'_\ka W_\K \int_1^z \left(\frac{1}{\alpha_\ka}-S_\ka(t)\right)\frac{P_\K\left(\log\frac{z}{t}\right)}{t}\;dt\notag\\
				 &= \frac{C'_\ka W_\K}{\alpha_\ka} \int_1^z \left(1-\alpha_\ka S_\ka(t)\right)\frac{P_\K\left(\log\frac{z}{t}\right)}{t}\;dt.
\end{align}
Note that $\frac{1}{\alpha_k}$ can be written as a Euler product in the following way
\begin{align*}
\frac{1}{\alpha_{\ka}} &=\sum_{\substack{d=1\\ (d,P(\ka))=1}}^\infty \mu^2(d)\frac{\tau_\K(d)}{\psi_\ka(d)d} = \prod_{p>\ka}\left(1 + \frac{\tau_\K(p)}{\psi_\ka(p)p}\right)\\
                          &=\prod_{p>\ka}\frac{p^2-(\ka)p+\K}{p(p-(\ka))} = \prod_{p>\ka}\left(1-\frac{1}{p}\right)\frac{p-\K}{p-(\ka)}.
\end{align*}
This allows us to express the leading constant as
\begin{align*}
&W_\K C'_\ka\frac{1}{\alpha_\ka}\\
=& \prod_p \left(1-\frac{\chi_\K(p)}{p}\right)^{-1}\left(1-\frac{1}{p}\right)^{\K}\prod_{p\leq\K}\left(1-\frac{1}{p}\right)\frac{P(\K)}{P(\ka)}\\
 &\prod_{p>\ka}\left(1-\frac{1}{p}\right)\frac{p-\K}{p-(\ka)}\prod_{p\leq \ka}\left(1-\frac{1}{p}\right)\\
=& \prod_p \left(1-\frac{\chi_\K(p)}{p}\right)^{-1}\left(1-\frac{1}{p}\right)^{\K}\\
 &\prod_{p\leq\ka}\left(1-\frac{\chi_\ka(p)}{p}\right)^{-1}\left(1-\frac{\chi_\K(p)}{p}\right)\left(1-\frac{1}{p}\right)\\
 &\prod_{p>\ka}\left(1-\frac{\chi_\ka(p)}{p}\right)^{-1}\left(1-\frac{\chi_\K(p)}{p}\right)\left(1-\frac{1}{p}\right)\prod_{p\leq \ka}\left(1-\frac{1}{p}\right)\\
=& \prod_p \left(1-\frac{\chi_\ka(p)}{p}\right)^{-1}\left(1-\frac{1}{p}\right)^{\ka}\prod_{p\leq\ka}\left(1-\frac{1}{p}\right) = W_\ka.
\end{align*}
So inequality (\ref{EQ3}) becomes
\begin{align}\label{EQ4}
F_\ka(z) &\geq W_\ka\left(Q_\K(\log z) -\int_1^z\alpha_\ka S_\ka(t)\frac{P_\K\left(\log\frac z t\right)}{t}\;dt \right)\notag\\
         &=W_\ka\left(Q_\K(\log z) -\int_1^z\alpha_\ka S_\ka(t)\left(\frac{\log^\K\frac{z}{t}}{\K !}-C_\K\log^{\K-1}\frac{z}{t}\right)\;\frac{dt}{t} \right)\notag\\
				 &\geq W_\ka\left(Q_\K(\log z) -\int_1^zS_\ka(t)\frac{\log^\K\frac{z}{t}}{\K !}\;\frac{dt}{t} \right)
\end{align}
Now we turn our attention to $S_k(t)$. The main difficulty in estimating this sum is caused by those $d$ who have small prime divisors in the range $k<p\leq k^2$. Therefore, we consider 
\[ \tilde{S}_k(z) = \sum_{\substack{d>z\\ (d,P(k^2))=1}}\mu^2(d)\frac{\tau_{k-1}(d)}{\psi_k(d)d}.\]
first. Clearly, each $d$ in this sum is only divided by primes $p>k^2$. For those $p$ we have 
\[\frac{k-1}{p-k} < \frac{k}{p}.\]
Employing Rankin's trick with some $0 < \lambda < 1$ yields
\begin{align}\label{EQ5}
\tilde{S}_k(z) &< \sum_{\substack{d>z\\ (d,P(k^2))=1}} \mu^2(d)\frac{\tau_k(d)}{d^2}\notag\\
               & = \frac{1}{z^{1-\lambda}}\sum_{\substack{1>z/d\\ (d,P(d^2))=1}}\mu^2(d)\frac{\tau_k(d)}{d^{1+\lambda}}\left(\frac{z}{d}\right)^{1-\lambda}\notag\\
							 & \leq\frac{1}{z^{1-\lambda}}\sum_{d:(d,P(k^2))=1}\mu^2(d)\frac{\tau_k(d)}{d^{1+\lambda}}\notag\\
							 &= \frac{1}{z^{1-\lambda}}\prod_{p>k^2}\left(1+\frac{k}{p^{1+\lambda}}\right)\leq\frac{1}{z^{1-\lambda}}\exp\left(\sum_{p>k^2}\frac{k}{p^{1+\lambda}}\right)\notag\\
							 &\leq\frac{1}{z^{1-\lambda}}\exp\left(k\int_{k^2}^\infty\frac{dt}{t^{1+\lambda}}\right) = \frac{1}{z^{1-\lambda}}\exp\left(\frac{1}{\lambda}k^{1-2\lambda}\right) 
\end{align}
for $z>0$. We will choose the exact value of the parameter $\lambda$ later. Next, we express $S_k(z)$ in terms of $\tilde{S}_k(z)$, which, combined with (\ref{EQ5}), gives us the following
\begin{align*}
S_k(z) &= \sum_{l\mid \frac{P(k^2)}{P(k)}} \mu^2(l)\frac{\tau_{k-1}(l)}{\psi_k(l)l}\tilde{S}_k\left(\frac{z}{l}\right)\\
       &\leq \sum_{l\mid \frac{P(k^2)}{P(k)}} \mu^2(l)\frac{\tau_{k-1}(l)}{\psi_k(l)l}\left(\frac{\ell}{z}\right)^{1-\lambda}\exp\left(\frac{1}{\lambda}k^{1-2\lambda}\right)\\
			 &= \frac{1}{z^{1-\lambda}}\exp\left(\frac{1}{\lambda}k^{1-2\lambda}\right)\sum_{l\mid \frac{P(k^2)}{P(k)}} \mu^2(l)\frac{\tau_{k-1}(l)}{\psi_k(l)l^\lambda}\\
			 &= \frac{1}{z^{1-\lambda}}\exp\left(\frac{1}{\lambda}k^{1-2\lambda}\right)\prod_{k<p\leq k^2}\left(1+\frac{k-1}{(p-k)p^\lambda}\right)\\
			 &\leq \frac{1}{z^{1-\lambda}}\exp\left(\frac{1}{\lambda}k^{1-2\lambda}\right)\prod_{k<p\leq k^2}\left(1+\frac{k^{1-\lambda}}{(p-k)}\right)\\
			 &\leq \frac{1}{z^{1-\lambda}}\exp\left(\frac{1}{\lambda}k^{1-2\lambda}\right)\exp\left(k^{1-\lambda}\left(1+\int_{k+1}^{k^2}\frac{dt}{(t-k)^2}\right)\right)\\
			 &\leq \frac{1}{z^{1-\lambda}}\exp\left(\frac{1}{\lambda}k^{1-2\lambda}\right)\exp\left(2k^{1-\lambda}\right).
\end{align*}
Using this estimate for (\ref{EQ4}) leads us to 
\begin{align*}
F_\ka(z)\geq W_\ka\left(Q_\K(\log z)-\frac{\exp\left(\frac{1}{\lambda}(\ka)^{1-2\lambda}\right)\exp\left(2(\ka)^{1-\lambda}\right)}{\K !}\right.\times\\
\left.\int_1^z\frac{\log^\K\frac{z}{t}}{t^{2-\lambda}}\;dt\right).
\end{align*}
The occurring integral is fairly easy to handle
\begin{align*}
\int_1^z\frac{\log^\K\frac{z}{t}}{t^{2-\lambda}}\;dt \leq \log^\K z\int_1^z\frac{dt}{t^{2-\lambda}}\leq \frac{\log^\K z}{1-\lambda}.
\end{align*}
Thus, together with the explicit expression of $Q_\K$, we have 
\begin{align*}
F_\ka(z)&\geq W_\ka\left(\frac{1}{(\ka)!}\log^\ka z-\frac{C_\K}{\K}\log^\K z\right.\\
&\quad\quad\left.-\frac{\exp\left(\frac{1}{\lambda}(\ka)^{1-2\lambda}\right)\exp\left(2(\ka)^{1-\lambda}\right)}{\K !}\frac{\log^\K z}{1-\lambda}\right).
\end{align*}
Finally, choosing $\lambda = 1-\frac{1}{\log\ka}$ allows us to write
\[ \frac{C_\K}{\K}+ \frac{\exp\left(\left(1-\frac{1}{\log\ka}\right)^{-1}\frac{e^2}{\ka}\right)\exp\left(2e\right)}{\K !}\log(\ka) = C_\ka\]
and the induction is completed.
\end{proof}
\subsection{A upper bound for \texorpdfstring{$S_3(z)$}{S3(z)}}

The estimation of $S_3$ follows similar steps as that of $S_k$, except that we avoid Rankin's trick and calculate some constants more carefully. As before, we deal with $\tilde{S}_3(z)$ first. If we define $T'_k(z)=\sum_{\substack{d\leq z\\(d,P(9))=1}}\tau_k(d)$, this sum is bounded by
\begin{align*}
\tilde{S}_3(z) \leq \sum_{\substack{d>z\\ (d,P(9))=1}}\mu^2(d)\frac{\tau_3(d)}{d^2}\leq \sum_{\substack{d>z\\ (d,P(9))=1}}\frac{\tau_3(d)}{d^2}=-\frac{1}{z^2}T'_{3}(z) + 2\int_z^\infty T'_{3}(t)\;\frac{dt}{t^3}.
\end{align*}
A simple calculation shows that
\[ \sum_{\substack{d\leq z\\ (d,P(9))=1}} 1 \leq \frac{\varphi(P(9))}{P(9)}\left(z + \frac{53}{8}\right)\]
holds for all $z\geq 1$. Using the general recursion
\begin{align*}
T'_{k+1}(z) &= \sum_{\substack{d\leq z\\ (d,P(9))=1}} \tau_k(d) = \sum_{\substack{d\leq z\\ (d,P(9))=1}} \tau_{k-1}(d) \sum_{\substack{l\leq z/d\\ (l,P(9))=1}} 1 \\
            &\leq \frac{\varphi(P(9))}{P(9)}\sum_{\substack{d\leq z\\ (d,P(9))=1}}\tau_{k-1}(d)\left(\frac{z}{d} + \frac{53}{8}\right)\\
            &= \frac{\varphi(P(9))}{P(9)}\left(\left(\frac{53}{8}+1\right) T'_{k}(z) + z \int_1^z T'_{k}(t)\;\frac{dt}{t^2}\right) \\
\end{align*}
we get 
\[T'_3(z)\leq z\left(\frac{\varphi(P(9))}{P(9)}\right)^3\left(\log z + \frac{175}{8}\right)^2 + \left(\frac{53}{8}\frac{\varphi(P(9))}{P(9)}\right)^3.\]
This allows us to estimate $\tilde{S}_3(z)$ by
\begin{align*}
\tilde{S}_3(z) \leq& -\frac{1}{z^2}T_{3,P(9)}(z) +  \left(\frac{\varphi(P(9))}{P(9)}\right)^3\int_z^\infty \frac{\left(\log t + \frac{175}{8}\right)^2}{t^2}\;dt \\
                   &+\frac{2}{z^2}\left(J_{P(9)}\frac{\varphi(P(9))}{P(9)}\right)^3 \\
							 \leq& \frac{1}{z^2}\left(6.945 - T_{3,P(9)}(z)\right) +\frac{512}{42875}\frac{\log^2 z + 47.75\log z + 524.266}{z} \\
							 \leq& \frac{512}{42875}\frac{\log^2 z + 47.75\log z + 524.266}{z}
\end{align*} 
for $z\geq 13$. In fact, one can show by direct computation that this is already true for $z\geq 1$ and since the right side is monotonically increasing as $z = \frac{1}{35}$ is approached, the inequality holds for $z\geq \frac{1}{35}$ as well. This gives us the desired result, since we have
\begin{align*}
S_3(z)    =& \sum_{\ell\mid 35}\mu^2(\ell)\frac{\tau(\ell)}{\ell\psi_3(\ell)}\tilde{S}_3\left(\frac{z}{\ell}\right) \\
       \leq& \sum_{\ell\mid 35}\mu^2(\ell)\frac{\tau(\ell)}{\psi_3(\ell)}\frac{512}{42875}\frac{\log^2\frac{z}{\ell} + 45.75\log\frac{z}{\ell} + 524.266}{z}\\
			 \leq& \frac{1534}{42875}\frac{\log^2 z + 45.75\log z + 524.266}{z}
\end{align*}
for $z\geq 1$.


\section{Bounds for \texorpdfstring{$L(z,\omega_0)$}{L(z,w0)} and \texorpdfstring{$G(z,\omega)$}{G(z,w)}}

With Lemma \ref{Lem3} in hand, we are now able to estimate $L(z,\omega_0)$ and finally $G(z,\omega)$. The starting point is the inequality (\ref{EQ2}). Again, we begin with the case $\K=3$. We have 
\begin{align*}
&L(z,\omega_0) \geq \prod_{p\leq 3}\frac{p}{p-\omega(p)} F_3\left(\frac{z}{6}\right)\\
              &\geq W_3\prod_{p\leq 3}\frac{p}{p-\omega(p)}\frac{1}{6}\left(\log^3 \frac{z}{6} + 6C_{1,3}\log^2\frac{z}{6} + 6C_{2,3}\log\frac{z}{6} + 6C_{3,3}\right)\\
							&=\frac{1}{6}\prod_p\left(1-\frac{\omega_0(p)}{p}\right)^{-1}\left(1-\frac{1}{p}\right)^3\left(\log^3 z + D_{1,3}\log^2 z + D_{2,3}\log z + D_{3,3}\right),
\end{align*}
for $z\geq 6$ with $D_{1,3} = -14.5859$, $D_{2,3} = 67.755$ and $D_{3,3} = 20.270$. In combination with Corollary \ref{Cor1} we get 
\[
L(z,\omega) \geq\frac{1}{6}\prod_p\left(1-\frac{\omega(p)}{p}\right)^{-1}\left(1-\frac{1}{p}\right)^3\left(\log^3 z + D_{1,3}\log^2 z + D_{2,3}\log z + D_{3,3}\right).
\]
Next, we are turning our attention towards $G(z,\omega)$. Summation by parts and the choice of $N = \frac{3}{2}z^2$ results in 
\begin{align*}
G(z,\omega) =& \sum_{q\leq z}\mu^2(q)\left(1+\frac{3}{2}zqN^{-1}\right)^{-1}\prod_{p\mid q}\frac{\omega(p)}{p-\omega(p)} \notag\\
            =& \frac{L(z,\omega)}{1+\frac{3}{2}z^2N^{-1}} + \int_1^z L(t,\omega)\frac{\frac{3}{2}zN^{-1}}{\left(1+\frac{3}{2}ztN^{-1}\right)^2} dt \notag\\
						=& \frac{L(z,\omega)}{2} + z\int_1^z \frac{L(t,\omega)}{(z+t)^2} dt \notag\\
				 \geq& \frac{1}{6}\prod_p\left(1-\frac{\omega(p)}{p}\right)^{-1}\left(1-\frac{1}{p}\right)^3\Biggl(\frac{1}{2}\left(\log^3 z + \sum_{\ell=1}^3 D_{\ell,3}\log^{3-\ell}z\right) +\notag\\
             &z\int_{z_0}^z\frac{\log^3 t + D_{1,3}\log^2 t + D_{2,3}\log t + D_{3,3}}{(z+t)^2} dt\Biggr)
\end{align*}
for $z\geq 6$ (or $N\geq 54$). Note that we have $z\geq z_0\geq 6$ (we will choose $z_0$ in a moment) to ensure that the integrand is positive. The integral itself can be broken down into integrals of the form $I_n(z,z_0) = z\int_{z_0}^z \frac{\log^n t}{(z+t)^2}\;dt$. While the cases $n=0,1$ are straightforward, we use  
\begin{align}\label{EQ6}
z\int_{z_0}^z\frac{\log^n t}{(z+t)^2} \;dt & = \frac{\log^n z}{2} - n\log 2\log^{n-1} z - \frac{z_0}{z+z_0}\log^n z_0 + \notag\\ 
                                           &n\log\frac{z+z_0}{z}\log^{n-1}z_0 + n(n-1)\int_{z+z_0}^{2z}\frac{\log\frac{t}{z}\log^{n-2}(t-z)}{t(t-z)} \;dt.
\end{align}
for $n\geq 2$. This identity is obtained by using integration by parts twice. We deduce 
\begin{align}\label{EQ7}
G(z,\omega)\geq&\frac{1}{6}\prod_p\left(1-\frac{\omega(p)}{p}\right)^{-1}\left(1-\frac{1}{p}\right)^3\times\notag\\
               &\left(\log^3z + E_{1,3}\log^2z + E_{2,3}\log z + E_{3,3}\right)
\end{align}
for $z\geq z_0\geq 6$ and with
\begin{align*}
E_{1,3} &= D_{1,3} - 3\log2\geq -16.666\\
E_{2,3} &=D_{2,3}-2D_{1,3}\log2 \geq 87.975\\
E_{3,3} &=D_{3,3} - D_{2,3}\log 2 - \frac{z_0}{z_0+z}\left(\log^3z_0 + D_{1,3}\log^2z_0 + D_{2,3}\log z_0 + D_{3,3}\right) \\
        &\quad+ R(z,z_0).
\end{align*}
Here the remainder in $E_{3,3}$ is given by $R(z,z_0) = \log\frac{z+z_0}{z}(3\log^2 z_0 + 2D_{1,3}\log z_0 + D_{2,3}) + 2\int_{z_0+z}^{2z}\frac{\log\frac{t}{z}}{t(t-z)}(3\log(t-z)+D_{1,3})\;dt$. Now we choose $z_0 = 130$ so that this integrand is also positive and thus negligible. If we additionally require $ z\geq z_1 = 540$, then $R(z,z_0) \geq -0.68187$ and we calculate $E_{3,3} = -50.624$. Inserting our choice $z=\sqrt{\frac{2}{3}N}$ into (\ref{EQ7}) proves ultimately the first part of Theorem \ref{THM2} for $N\geq 437504\geq \frac{3}{2}z_1^2$. To show the remaining part, we turn back to inequality (\ref{EQ2}) and lemma \ref{Lem3}. It follows for $\K\geq 3$ and $z> P(\K)$ that
\begin{align*}
L(z,\omega_0) &=\sum_{\ell\mid P(\K)}\mu^2(\ell)\prod_{p\mid \ell}\frac{\omega_0(p)}{p-\omega_0(p)} F_\K\left(\frac{z}{\ell}\right)\geq \prod_{p\leq\K}\left(1-\frac{\omega_0(p)}{p}\right)^{-1} F_\K\left(\frac{z}{P(\K)}\right)\\
              &\geq \prod_{p\leq\K}\left(1-\frac{\omega_0(p)}{p}\right)^{-1}\prod_p \left(1-\frac{\chi_\K(p)}{p}\right)^{-1}\left(1-\frac{1}{p}\right)^{\K}\prod_{p\leq\K}\left(1-\frac{1}{p}\right)\\
              &\quad\left(\frac{1}{\K!}\log^\K \frac{z}{P(k)} - C_\K\log^{\K-1} \frac{z}{P(k)}\right)\\
              &=\frac{1}{\K!}\prod_p\left(1-\frac{\omega_0(p)}{p}\right)^{-1}\left(1-\frac{1}{p}\right)^\K\log^\K z\\
              &\quad\left(\left(1-\frac{\log P(\K)}{\log z}\right)^\K-\frac{C_\K\K!}{\log z}\left(1-\frac{\log P(\K)}{\log z}\right)^{\K-1}\right)\\
              &\geq \frac{1}{\K!}\prod_p\left(1-\frac{\omega_0(p)}{p}\right)^{-1}\left(1-\frac{1}{p}\right)^\K\log^\K z\left(1-\frac{\K\vartheta(\K)}{\log z}-\frac{C_\K\K!}{\log z}\right)\\
\end{align*}
Of course, $\vartheta(x)$ is the well known first Chebyshev function. According to Theorem 9 from \cite{ROS}, it is bounded by
\begin{align}\label{EQ8}
\vartheta(x) = \sum_{p\leq x}\log p < 1.01624x.
\end{align}
for $x>1$. In combination with Corollary \ref{Cor1} this gives us 
\begin{align*}
L(z,\omega)&\geq \prod_{\substack{p\in\mathfrak{P}_0\\ p>\K}}\frac{p-\omega_0(p)}{p-\omega(p)}L(z,\omega_0)\\
           &\geq \frac{1}{\K!}\prod_p\left(1-\frac{\omega(p)}{p}\right)^{-1}\left(1-\frac{1}{p}\right)^\K\log^\K z\left(1-\frac{C'_\K}{\log z}\right),
\end{align*} 
with $C'_\K = 1.02\K^2 + \K!C_\K$. Next we deal with $G(z,\omega)$ again. We will follow the same lines as in the proof of the special case $\K=3$. Specifically, we will choose $z=\sqrt{\frac{2}{3}N}$ and use summation by parts. Thus, we have 
\begin{align*}
G(z,\omega) &= \frac{L(z,\omega)}{2} + z\int_1^z \frac{L(t,\omega)}{(t+z)^2}\; dt\\
            &\geq \frac{1}{\K!}\prod_p\left(1-\frac{\omega(p)}{p}\right)^{-1}\left(1-\frac{1}{p}\right)^\K\times\\
						&\quad\quad\left(\frac{\log^\K z}{2}\left(1-\frac{C'_\K}{\log z}\right) + z\int_{\exp(C'_\K)}^z \frac{\log^\K t}{(t+z)^2}\left(1-\frac{C'_\K}{\log t}\right)\; dt\right)\\
						&\geq \frac{1}{\K!}\prod_p\left(1-\frac{\omega(p)}{p}\right)^{-1}\left(1-\frac{1}{p}\right)^\K\times\\
						&\quad\quad\left(\frac{\log^\K z}{2}\left(1-\frac{C'_\K}{\log z}\right) + z\int_{\exp(C'_\K)}^z \frac{\log^\K t}{(t+z)^2}\; dt-C'_\K\frac{\log^{\K-1}z}{2} \right)\\
\end{align*}
for $z\geq \exp(C'_\K)$. We derive the lower bound
\begin{align*}
z\int_{\exp(C'_\K)}^z \frac{\log^\K t}{(t+z)^2} dt &\geq \frac{\log^\K z}{2} - \K\log2\log^{\K-1}z - \frac{e^{C'_\K}}{z+e^{C'_\K}}(C'_\K)^\K\\
                                                   &\geq \frac{\log^\K z}{2} - \left(\frac{1}{2}C'_\K + \K\log2\right)\log^{\K-1}z
\end{align*}
for $z\geq\exp{C'_\K}$ and $\K\geq 3$ using the identity (\ref{EQ6}) from our first encounter with integrals of the from $I_n(z,z_0)$. This gives us
\begin{align*}
G(z,\omega)\geq \frac{1}{\K!}\prod_p\left(1-\frac{\omega(p)}{p}\right)^{-1}\left(1-\frac{1}{p}\right)^\K\left(\log^\K z - \left(\K\log2 + \frac{3}{2}C'_\K\right)\log^{\K-1} z\right).
\end{align*}
If we insert our choice of $z = \sqrt{\frac{2}{3}N}$ into this expression, we derive the following inequality in terms of $N$
\begin{align*}
G(z,\omega) \geq& \frac{1}{\K!2^\K}\prod_p\left(1-\frac{\omega(p)}{p}\right)^{-1}\left(1-\frac{1}{p}\right)^\K\log^\K N\times\\
                &\left(\left(1+\frac{\log\frac{2}{3}}{\log N}\right)^\K - \frac{2\left(\K\log2+\frac{3}{2} C'_\K\right)}{\log N}\right)\\
						\geq& \frac{1}{\K!2^\K}\prod_p\left(1-\frac{\omega(p)}{p}\right)^{-1}\left(1-\frac{1}{p}\right)^\K\log^\K N\times\\
						    &\left(1+\frac{\K\log\frac{2}{3}-2(\K\log2+3 C'_\K)}{\log N}\right).
\end{align*}
Using Sterling's well known approximation formula 
\[
\log k!\leq k\log k - k +\frac 1 2 \log k+ \frac 1 2 \log 2\pi +\frac{1}{12k}
\]
for $k\geq 1$, we can bound the constant $R_\K$ of the error term by
\begin{align*}
R_\K &= \K\log 6 + 3\left(1.02\K^2+\K !C_\K\right)\\
     &\leq \K\log 6 +3\left(1.02\K^2 + \K\exp\left(2e+\frac{e^2}{3}\left(\frac{\log 3}{\log 3-1}\right)\right)\log\K !\right)\\
		 &\leq 3\left(1.02\K^2 + \K^2\log\K\exp\left(2e+\frac{e^2}{3}\left(\frac{\log 3}{\log 3-1}\right)\right)\right)\\
		 &\leq C\K^2\log\K = r_\K,
\end{align*}
with $C = 5.961\cdot 10^{14}$. The function $G(z,\omega)$ is therefore bounded through
\[
G(z,\omega)\geq \frac{1}{\K!2^\K}\prod_p\left(1-\frac{\omega(p)}{p}\right)^{-1}\left(1-\frac{1}{p}\right)^\K\log^\K N\left(1 - \frac{r_\K}{\log N}\right)
\]
for $\K\geq 3$ and with $\sqrt{\frac{2}{3}N} = z\geq \exp(C'_\K)$. But the right side is positive only for $N>\exp(r_\K)\geq\exp(\K\log 6 + 3 C'_\K) \geq \frac{3}{2}\exp(2C'_\K)$. Thus, lemma \ref{Lem1} is applicable and we derive the desired result.


\section{Application of the \texorpdfstring{$K$}{K}-dimensional sieve}

First we prove the general case for $K\geq 3$ and all admissible configurations $\bs{a},\bs{b}\in\mathbb{Z}^{K-1}$. For $z\geq 2$ we have
\begin{align}\label{EQ9}
\pi_{K}(\bs{a},\bs{b};X)\leq \pi(z)+\sum_{\substack{z<m\leq X\\ (m\prod(a_im+b_i),P(z))=1}}1.
\end{align}
We see that the second sum can be expressed as a sieving function. In fact, let 
\begin{align*}
\mathcal{A}=\{z < m\leq x\mid m\in\mathbb{N}\}
\end{align*}
and $\Omega_p = \Omega_p(\bs{a},\bs{b})$ be the set of all solutions $\pmod{p}$ of the congruence
\begin{align}\label{EQ10}
m\prod_{i=1}^{K-1}(a_im+b_i) \equiv 0\pmod{p}.
\end{align}
Then we can write the above sum as 
\begin{align*}
\sum_{\substack{z<m\leq X\\ (m\prod(a_im+b_i),P(z))=1}}1 &= \#\{z<m\leq x\mid m\pmod{p}\not\in \Omega_p(\bs{a},\bs{b}),\;\forall p\leq z\}\\
                                                         &= S_K(\mathcal{A},\Omega,z)
\end{align*}
Thus, $\rho_K(p) = \vert\Omega_p\vert$ is in fact the number of solutions of (\ref{EQ10}) and we have
\begin{align*}
\rho_{K}(p) &= K\quad\text{for}\; p\nmid E_{K}\\
\rho_{K}(p) &< K\quad\text{for}\; p\mid E_{K}.
\end{align*}
This means that $\rho_K(p)$ satisfies (\ref{VOR1}) with the exception-set $\Ausn = \{p\in\mathbb{P}\mid p\mid E_{K+1}\}$. Note that we have $\rho_{K}>0$, since every $m\equiv 0\;(\bmod{p})$ is a trivial solution of (\ref{EQ10}).  Also, $\rho_{K}(p)$ fulfils the condition (\ref{VOR2}) due to the requirements of our Theorem. Therefore, Theorem \ref{THM2} is applicable and we have 
\begin{align}\label{EQ11}
S_K(\mathcal{A},\Omega,z) <2^{K}K!&\prod_p\left(1-\frac{\rho_{K}(p)}{p}\right)\left(1-\frac{1}{p}\right)^{-K}\frac{X}{\log^{K} X}\left(\frac{\log X}{\log X-r_{K}}\right)
\end{align}
for $z=\sqrt{\frac{2}{3}X}$ and $X>\exp(r_{K})$. For convenience we define the product over primes as
\[\tilde{W}_{K} = \prod_p\left(1-\frac{\rho_{K}(p)}{p}\right)\left(1-\frac{1}{p}\right)^{-K}.\]
Using the estimate (\ref{EQ11}) and the prime number Theorem 
\[\pi(x) < 1.25506 \frac{x}{\log x},\quad(x>1)\]
(see \cite{ROS}, Corollary 1), we have
\begin{align}\label{EQ12}
\pi_{K}(\bs{a},\bs{b};X) &\leq \pi(z) + S_K(\mathcal{A},\Omega,z)\notag\\
                         &\leq 2^KK!\tilde{W}_K\left(\frac{1.3}{2^KK!\tilde{W}_K}\frac{\sqrt{\frac{2}{3}X}}{\frac{1}{2}\log\frac{2}{3}X} + \frac{X}{\log^KX}\frac{\log X}{\log X - r_K}\right).
\end{align}
Next, we require suitable bounds for the product $\tilde{W}_K$, which we obtain by the following 
\begin{lem}\label{Lem5}
Let $K\geq 2$ and $\bs{a}_{K-1},\bs{b}_{K-1}\in\mathbb{Z}^{K-1}$ be admissible in the sense of Theorem \ref{THM1}. We define $\bs{a}_{K-2},\bs{b}_{K-2}\in\mathbb{Z}^{K-2}$ through
$(a_i)_{K-2} = (a_i)_{K-1}$, $(b_i)_{K-2} = (b_i)_{K-1}$ for $1\leq i < K-1$. Then we have 
\[\frac{\tilde{W}_K}{\tilde{W}_{K-1}} >\frac{1}{2K}e^{-1.3K}\]
for $\rho_{i+1}(p) = |\Omega_p(\bs{a}_i,\bs{b}_i)|$.
\end{lem}
Note that the vectors $\bs{a}_{K-2},\bs{b}_{K-2}$ are of course also admissible, but for the sieve dimension $K-1$. Thus, $\rho_{K-1}$ has the same properties as $\rho_K$.
\begin{proof}
We start by writing the quotient as
\begin{align}\label{EQT4}
\frac{\tilde{W}_K}{\tilde{W}_{K-1}} &= \prod_p\frac{p-\rho_K(p)}{p-\rho_{K-1}(p)}\left(1-\frac{1}{p}\right)^{-1}\notag\\
          &=\prod_{p\leq K}\frac{p-\rho_K(p)}{p-\rho_{K-1}(p)}\frac{p}{p-1}\prod_{p>K}\frac{p-\rho_K(p)}{p-\rho_{K-1}(p)}\frac{p}{p-1}\notag\\
          &=\mathcal{P}_1\cdot\mathcal{P}_2.
\end{align}
Because of $1\leq \rho_{K-1}(p)\leq \rho_K(p)\leq p-1$ we have
\begin{align*}
\mathcal{P}_1 &= \prod_{p\leq K}\frac{p-\rho_K(p)}{p-\rho_{K-1}(p)}\frac{p}{p-1} \geq \prod_{p\leq K}\frac{p}{(p-1)^2}\\
              &= \exp\left(\sum_{p\leq K}\log p-2\log(p-1)\right)\\
              &= \exp\left(\pi(K)\log\frac{K}{(K-1)^2} -\pi(2)\log2- \int_2^{K}\pi(t)\left(\frac{1}{t}-\frac{2}{t-1}\right)\;dt\right)\\
              &= \exp\left(\pi(K)\log\frac{K}{(K-1)^2}-\log2+ \int_2^{K}\pi(t)\frac{t+1}{t(t-1)}\;dt\right)\\
              &\geq \exp\left(\pi(K)\log\frac{K}{(K-1)^2} - \log2\right).
\end{align*}
With the help of the prime number Theorem from earlier, we get
\[\mathcal{P}_1 > \exp\left(\pi(K)\log\frac{K}{(K-1)^2}-\log2\right) > \exp(-5.5K)\]
for $K\geq 3$. However, the inequality holds for $K=2$ as well. The second factor $\mathcal{P}_2$ can be estimated by
\begin{align*}
\mathcal{P}_2 &= \prod_{p>K}\frac{p-\rho_{K}(p)}{p-\rho_{K-1}(p)}\frac{p}{p-1}\\
              &= \prod_{\substack{p>K\\ \rho_K(p)=\rho_{K-1}(p)+1}}\frac{p-\rho_K(p)}{p-\rho_{K-1}(p)}\frac{p}{p-1}\prod_{\substack{p>K\\ \rho_K(p)=\rho_{K-1}(p)}}\frac{p-\rho_K(p)}{p-\rho_{K-1}(p)}\frac{p}{p-1}\\
              &> \prod_{\substack{p>K \\ \rho_K=\rho_{K-1}+1}}\frac{p-K}{p-(K-1)}\frac{p}{p-1}> \prod_{p>K}\frac{p-K}{p-(K-1)}\frac{p}{p-1}\\
              &\geq \exp\left(\sum_{d>K}\log\frac{d-K}{d-(K-1)}+\log\frac{d}{d-1}\right)\\
              &\geq \exp\left(\int_{K+1}^\infty\log\frac{t-K}{t-(K-1)}+\log\frac{t}{t-1}\;dt\right)\\
              &=\exp\Bigg(2\log2-2-\log\left(1-\frac{1}{K+1}\right)^{K+1}-\log K\Bigg)\\
              &>\exp\left(2\log2-1-\log K\right) > \frac{1}{K}>\frac{1}{2K}.
\end{align*}
The desired result is obtained by combining both estimates.
\end{proof}
As a simple consequence we have
\[
\tilde{W}_K = \prod_{i=2}^K \frac{\tilde{W}_i}{\tilde{W}_{i-1}} >\frac{1}{2^{K-1}K!}e^{-5.5K^2}.
\]
for $K\geq 2$. This yields for (\ref{EQ12})
\begin{align*}
\pi_{K}(\bs{a},\bs{b};X) \leq& 2^KK!\tilde{W}_K\frac{X}{\log^KX}\left(1.3\sqrt{\frac{2}{3}}\left(\frac{\log X}{\log\frac{2}{3}X}\right)e^{5.5K^2}\frac{\log^{K-1}}{\sqrt{X}} + \frac{\log X}{\log X - r_K}\right)\\
                            <& 2^KK!\tilde{W}_K\frac{X}{\log^KX}\left(1.08e^{5.5K^2}\frac{\log^{K-1}}{\sqrt{X}} + \frac{\log X}{\log X - r_K}\right)\\
												    =& 2^KK!\tilde{W}_K\frac{X}{\log^KX}\left(\frac{\log X}{\log X - 2r_K}\right)\times\\
														 &\left(1.08e^{5.5K^2}\frac{\log^{K-2}}{\sqrt{X}}\left(\log X - 2r_K\right) +  1 - \frac{r_K}{\log X - r_K}\right)
\end{align*}
for $X>\exp(2r_K)$. If we can show that
\[1.08e^{5.5K^2}\frac{\log^{K-1}X}{\sqrt{X}} < \frac{r_K}{\log X - r_K}\]
is true for $X>\exp(2r_K)$, then the first part of Theorem \ref{THM1} follows. First, we have
\[\log X < X^{\frac{1}{4(K-1)}}\]
for such $X$ and $K\geq 3$. So the left side is bounded from above by 
\[1.08e^{5.5K^2}\frac{\log^{K-1}X}{\sqrt{X}} < 1.08e^{5.5K^2}\frac{1}{X^{\frac 1 4}}.\]
But $X^{1/4}$ grows at a faster rate than $\log X$ as $X\rightarrow \infty$. Therefore, it is sufficient to check whether 
\[ 1.08 e^{5.5K^2}\frac{1}{e^{\frac{1}{2}r_K}} < \frac{r_K}{\log e^{2r_K} -r_K} = 1\]
is true. In fact we have $1.08e^{5.5K^2-0.5r_K} = 1.08e^{5.5K^2 - 0.5CK^2\log K} < 1/2$, since $C$ is much larger than $5.5$. Thus, the case $K\geq 3$ follows immediately. For $K=3$, we start by writing
\begin{align*}
\pi_3(\bs{a},\bs{b};X) &\leq \pi_2(a_1,b_1;z) + \sum_{(m\prod(a_im+b_i),P(z))=1} 1\\
                       &= \pi_2(a_1,b_1;z) + S_3(\mathcal{A},\Omega,z)
\end{align*}
for $z\geq 2$ and admissible $\bs{a},\bs{b}\in\mathbb{Z}^2$. This time, we use the sharper estimate for $S_3$ from Theorem \ref{THM2} and the main result of Siebert
\[\pi_2(a,b,x) < 8\prod_p\left(1-\frac{\rho_2(p)}{p}\right)\left(1-\frac{1}{p}\right)^{-2}\frac{x}{\log^2x},\quad (x>1).\]
Choosing $z=\sqrt{\frac{2}{3}X}\geq 540$ gives us
\begin{align*}
\pi_3(\bs{a},\bs{b};X) \leq& 8\prod_p\left(1-\frac{\rho_2(p)}{p}\right)\left(1-\frac{1}{p}\right)^{-2}\frac{\sqrt{\frac{2}{3}X}}{\frac{1}{4}\log^2\frac{2}{3}X} + \\
                           & 48\prod_p\left(1-\frac{\rho_3(p)}{p}\right)\left(1-\frac{1}{p}\right)^{-3}\frac{X}{P(\log X)}\\
													=& 48\tilde{W}_3\frac{X}{\log^3X}\left(\frac{32}{48}\frac{\tilde{W}_2}{\tilde{W}_3}\sqrt{\frac{2}{3}}\frac{\log X}{\sqrt{X}}\left(\frac{\log X}{\log\frac{2}{3}X}\right)^2 + \frac{\log^3X}{P(\log X)}\right)\\
													\leq& 48\tilde{W}_3\frac{X}{\log^3X}\left(0.57998\frac{\tilde{W}_2}{\tilde{W}_3}\frac{\log X}{\sqrt{X}} + \frac{\log^3X}{P(\log X)}\right)
\end{align*}
for $X\geq 437400$. According to the proof of Lemma \ref{Lem5}, we have
\[\frac{\tilde{W}_2}{\tilde{W}_3} = \mathcal{P}_1\cdot\mathcal{P}_2 > \prod_{p\leq 3}\frac{p}{(p-1)^2}\cdot \frac{1}{3} = \frac 1 2.\]
It follows that
\[ \pi_3(\bs{a},\bs{b};X)\leq 48\tilde{W}_3\frac{X}{\log^3X}\left(1.15996\frac{\log X}{\sqrt{X}} + \frac{\log^3X}{P(\log X)}\right).\]
Our Theorem is shown, if the term in brackets is $\leq 1 + 34.549\log^2X/P(\log X)$ for $X\geq X_0 = 437400$. But this inequality is equivalent to 
\[1.15996\frac{\log^4X}{\sqrt{X}}\frac{P(\log X)}{\log^3{X}}\leq P(\log X) - \log^3X + 34.594\log^2X.\]
We see, that the left side is decreasing for $X\geq e^8 = 2980.957\ldots$, while the right side increases for $X>0$. So we only need to check if the inequality holds for $X=X_0$. Since this is indeed true, we have successfully completed the proof of Theorem \ref{THM1}.


\end{document}